\documentclass[12pt]{article}
\usepackage[utf8]{inputenc}
\usepackage[english]{babel}
\usepackage[top=1.5cm, outer=1.5cm, bottom=1.5cm, inner=1.5cm]{geometry}
\usepackage{amssymb}
\usepackage{amsmath}
\usepackage{amsthm}

\DeclareMathOperator{\IN}{in}

\DeclareMathOperator{\Char}{char}

\begin{document}

\begin{center}\Large{Prime differential nilalgebras exist.}\end{center}

\begin{center}{Gleb Pogudin\footnote{Moscow State University, {e-mail: \texttt{pogudin.gleb@gmail.com.}}}}\end{center}

\begin{abstract}
We construct a monomorphism from the differential algebra $k\{x\} / [x^m]$ to a Grassmann algebra endowed with a structure of a differential algebra. Using this monomorphism, we prove primality of $k\{x\} / [x^m]$ and its algebra of differential polynomials, solve one of Ritt's problems and give a new proof of integrality of the ideal $[x^m]$.


\end{abstract}

From now on, we will always assume that $k$ is a field of characteristic zero.

{\bf 1. The main construction. } The differential ideal $[x^m]$, $m\in\mathbb{N}$, first appeared in Levi's paper \cite{Levi}. He used combinatorial properties of $[x^m]$ in the proof of low power theorem. Levi proved a sufficient condition for a monomial in $k\{x\}$ to lie in $[x^m]$. Ritt conjectured (\cite{Ritt}) that a minimal $j$, such that $x_i^j\in [x^m]$, is given by a formula $(i + 1)m - i$. This conjecture was proved for $i \leqslant 2$ by O'Keefe \cite{Keefe}. We will prove this statement for all $i$ (Theorem 3).

As mentioned in \cite{Zobnin}, Levi's reduction process used in \cite{Levi} is similar to the reduction with respect to a differential Groebner basis.

We denote a differential polynomial algebra in one indeterminate by $k\{ x \}$. More precisely, $k\{x\}$ is a polynomial algebra $k[x_0, x_1, x_2, \ldots]$ in a sequence of algebraic independent indeterminates $(x = x_0, x_1, x_2, \ldots)$ endowed with a derivation (i.e. $k$-linear operator satisfying the Leibniz law) such that $x_n' = x_{n + 1}$. An ideal $I$ is called a differential ideal if $I'\subset I$. Let $k_{+}\{x\}$ be the subalgebra of the differential polynomials with zero constant term.

The main construction is introduced as follows. Let us consider an infinite dimensional vector space $V_m$ with a basis consisting of $\xi_i^k$ and $\eta_i^k$ ($k = 0, \dots, m - 2$, $i\in \mathbb{Z}_{\geqslant 0}$). We denote a Grassmann algebra of $V_m$ by $\Lambda(V_m)$ and denote its even and odd components by $\Lambda_0(V_m)$ and $\Lambda_1(V_m)$, respectively. Let us note that, from now on, $\Lambda(V_m)$ is assumed not to contain unity. Let us equip $\Lambda(V_m)$ with a derivation (not superderivation) using the formulas $(\xi_i^k)' = \xi_{i + 1}^k$ and $(\eta_i^k)' = \eta_{i + 1}^k$, i.e. the derivation increments the subscript. Obviously, $\Lambda_0(V_m)$ is a differential subalgebra in $\Lambda(V_m)$.

Let us denote the quotient algebra $k_{+}\{x\} / [x^m]$ by $D_m$ and the image of $x$ in $D_m$ by $\bar{x}$. We define a homomorphism of differential algebras $\varphi_m\colon D_m\rightarrow \Lambda_0(V_m)$ by the rule $\varphi_m(\bar{x}) = \sum\limits_{k = 0}^{m - 2}\xi_0^k\land\eta_0^k$.

{\bfseries Lemma 1. }\emph{$\varphi_m$ is a monomorphism.}

\begin{proof} Let us prove that $\varphi_m$ is a homomorphism. It is sufficient to prove that $(\varphi_m(\bar{x}))^m = 0$. Indeed, $(\varphi_m(\bar{x}))^m$ is a homogeneous skew-symmetric form of degree $2m$ in $2(m-1)$ variables.

    Let us recall that (according to Levi \cite{Levi}) a monomial $x^{p_0}x_1^{p_1}\cdots x_n^{p_n}$ ($p_i\in\mathbb{Z}_{\geqslant 0}$) in $k_{+}\{x\}$ is said to be an $\alpha_m$-monomial if $p_{i} + p_{i + 1} < m$ for each $i$. We need the following fact (\cite{Levi}, theorem 1.1):

{\bf Fact. } \emph{Images of $\alpha_m$-monomials in $D_m$ constitute a basis for the linear vector space $D_m$. }

Thus it remains to check that the kernel of $\varphi_m$ does not contain any non-trivial linear combination of $\alpha_m$-monomials. Let us introduce a monomial order. We fix two monomials $P = x^{p_0}\cdots x_n^{p_n}$ and $Q = x_0^{q_0}\cdots x_n^{q_n}$ ($p_i, q_i\in\mathbb{N}\cup\{0\}$). Let $j$ be the minimal index such that $p_j \neq q_j$. Then $P$ is said to be larger than $Q$ if $p_j < q_j$ and $Q$ is said to be larger than $P$ otherwise. This order is called $degrevlex$ order (\cite{Zobnin}).

Let us consider an arbitrary linear combination $L$ of $\alpha_m$-monomials and its leading term $M = \bar{x}_{k_n}\cdots \bar{x}_{k_0}$, where $k_n \geqslant\cdots\geqslant k_0$. The $\alpha_m$-property can be reformulated as follows: for each pair $i$ and $i'$ such that $i - i' = m - 1$ holds the inequality $k_i - k_{i'} > 1$. In particular, $k_i \geqslant 2\left[\frac{i}{m - 1}\right]$. We are going to find a monomial $\mu$ in $\varphi_m(M)$ which does not appear in $\varphi_m(M')$ for any monomial $M' \neq M$ from $L$. We assign to each $\bar{x}_{k_j}$ a monomial $\xi_{k_j - q}^r\land\eta_{q}^r$, where $q$ and $r$ are the quotient and the remainder of division $j$ by $m - 1$, respectively. Let $\mu$ be the product of these monomials. Due to the inequality above, $k_j - q \geqslant 0$. Moreover, the $\alpha_m$-property implies:

\begin{equation}
\label{chain}
k_j - q < k_{j + m - 1} - q - 1 < k_{j + 2m - 2} - q - 2 < \ldots
\end{equation}
for each $j$. Hence, all basis vectors appearing in $\mu$ are different and $\mu \neq 0$.

Let us prove that $\mu$ does not appear in any monomial $M' = \bar{x}_{k'_m}\cdots \bar{x}_{k'_0} < M$ ($k'_m\geqslant\cdots\geqslant k'_0$). 
Since $M' < M$, $k'_0 \geqslant k_0$. 
On the other hand, $\mu$ contains $\xi_{a}^{b}\land\eta_{c}^{b}$ from $\varphi_m(\bar{x}_{k'_0})$ for some $a$, $b$ and $c$. 
Inequality (\ref{chain}) implies that minimum possible value for $a + c$ is $k_0$, where $a = k_0$, $b$ equals the remainder of division $k_0$ by $m - 1$, and $c = 0$. 
Hence, $k_0 = k'_0$, $\varphi_m(\bar{x}_{k_0})$ and $\varphi_m(\bar{x}_{k'_0})$ correspond to the same submonomial of $\mu$. 
Continuing this line of reasoning, we see that $M = M'$.
Thus, $\varphi_m$ is injective.\end{proof}


{\bfseries 2. The algebra $k_{+}\{x\} / [x^2]$ is prime.}

{\bfseries Lemma 2. } \emph{Let $P$ be a prime differential algebra and $Q \subset P$ be its differential subalgebra. Then, $Q$ is prime.}

\begin{proof}
    Let $A$ and $B$ be differential ideals in $Q$ such that $AB = 0$. 
    Consider ideals (not necessary differential) $\bar{A}$ and $\bar{B}$ in $P$ generated by $A$ and $B$ respectively. 
    Obviously, $\bar{A}\bar{B} = 0$.
    On the other hand, we claim that $\bar{A}$ and $\bar{B}$ are actually differential ideals. 
    Indeed, let $a\in A$ and $x\in P$, then $(pa)' = p'a + pa'$, where $a, a'\in A$ and $p, p'\in P$. 
    This contradicts the primality of $P$.
\end{proof}


{\bfseries Theorem 1.} \emph{The algebra $D_2$ is prime. }

\begin{proof}

{\bfseries }
    Due to Lemma 2, it remains to check that $\Lambda_0(V_2)$ is prime. 
    Let $A$ and $B$ be such differential ideals that $AB = 0$. 
    Let us fix $a\in A$ and $b\in B$. 
    There exist monomials $m_a$ and $m_b$ such that $m_aa$ and $m_bb$ are nonzero monomials. 
    Hence, we can assume without loss of generality that $a$ and $b$ are monomials. 
    There also exists large enough $N$ such that $a^{(N)}$ contains a monomial with derivations of higher order than the highest order of the derivations in $b$. 
    Clearly, $a^{(N)}b \neq 0$.
\end{proof}

The above theorem gives us the example mentioned in the title of this paper. 
Indeed, all elements in $D_2$ are nilpotent because $\varphi_2$ imbeds $D_2$ into a Grassmann algebra. 
The example of an associative prime nilalgebra is presented in Theorem 1 of \cite{RybAnd}.

{\bfseries Corollary. }\emph{For any $a, b\in D_2$ there exists $k\in\mathbb{N}$ such that $b^{(k)}a\neq 0$.}

{\bfseries 3. The algebra of differential operators over $k_{+}\{x\} / [x^2]$ is prime.}

In paper \cite{Japan}, it was proved that the algebra of differential operators over a prime differential algebra with unity is prime. 

But $D_2$ does not contain the unity. 
By $(D_2)_{id}$, we denote $D_2$ with an externally adjoined unity. 
Let us prove the following technical lemma.

{\bfseries Lemma 3.} \emph{ Let $f(t_1, \ldots t_k)$ be a nonzero differential polynomial in variables $t_i$ ($1\leqslant i\leqslant k$) with coefficients in $(D_2)_{id}$. Then there exist $a_1,\ldots,a_k \in D_2$ such that $f(a_1,\ldots, a_k) \neq 0$.}

\begin{proof}
    Assume the converse. 
    Then, $f(t_1, \ldots, t_k)$ is a differential polynomial identity. 
    Since $\Char k = 0$, after a complete linearization of $f(t_1, \ldots, t_n)$ we obtain nonzero identity $\tilde{f}(\tilde{t}_1, \ldots, \tilde{t}_m)$ linear with respect to each of $\tilde{t}_i$.
    Moreover, without loss of generality it can be assumed that $m = 1$. 
    The general case follows from the case $m = 1$ by simple induction. 
    
    Thus, $f(t)$ is a linear form in $t$ and its derivations with coefficients in $(D_2)_{id}$ and $f(1) \neq 0$. 
    Let us represent the coefficients of $f$ as $\alpha_2$-polynomials. 
    By $N$, we denote the maximal order of a derivation of $x$ appearing in the coefficients of $f$ represented in such a way.  
    Let $a = \bar{x}_{N + 2}$. 
    Then $f(a)$ is a nontrivial linear combination of $\alpha_2$-monomials. 
    Hence $f(a) \neq 0$.
\end{proof}

By $A[\partial]$, we denote the algebra of differential operators over an algebra $A$.

{\bfseries Theorem 2.}\emph{The algebra $D_2[\partial]$ is prime.}

\begin{proof}

First of all, let us check that $(D_2)_{id}[\partial]$ is prime. 
It can be deduced from results of Paper \cite{Japan}, but, for the sake of completeness, we will provide a self-contained proof. 
Assume the converse, then there are ideals $A$ and $B$ in $(D_2)_{id}[\partial]$ such that $AB = 0$. 
Any $a \in (D_2)_{id}[\partial]$ is of the form $a = a_n\partial^n + \cdots + a_0$, where $a_i\in(D_2)_{id}$. 
By $\IN(A)$, we denote the set of all leading terms of an ideal $A$. 
Clearly, it is an ideal in $(D_2)_{id}$. 
Moreover, since $a\partial - \partial a = a'_n\partial^n + \cdots + a'_0$, it is a differential ideal. 
Hence $\IN(A)\IN(B) = 0$, where $\IN(A)$ and $\IN(B)$ are nonzero differential ideals in $(D_2)_{id}$. 
This contradicts  primality of $(D_2)_{id}$.

Assume that $D_2[\partial]$ is not prime. 
Hence, there exist $a, b \in D_2[\partial]$ such that $\langle a \rangle \langle b \rangle = 0$ (by $\langle x \rangle$, we denote the ideal generated by $x$), where $a = a_n\partial^n + \cdots + a_0$ and $b = b_m\partial^m + \cdots + b_0$. 
By the corollary above, there exists $k$ such that $(a_n)^{(k)}b_m \neq 0$. 
Let us consider $f(t) = [(t\partial)^k, a]b$. 
The leading term of $f(1)$ is $(a_n)^{(k)}b_m \neq 0$. 
Hence, by Lemma 3, there exists $c\in D_2$ such that $0 \neq [(c\partial)^k, a]b \in \langle a \rangle \langle b \rangle$.

\end{proof}

{\bfseries Proposition. }\emph{All elements of $D_2[\partial]$ are nilpotents.}

\begin{proof}
For each differential monomial $M = x_{k_n}\cdots x_{k_0}$, we assign the weight $w(M) = k_n + \cdots + k_0$. 
For each polynomial $p\in k\{x\}$, we denote the maximal weight of its terms by $w(p)$.
The weight of a monomial in $D_2$ can be defined in the same way. 
Obviously, the minimal weight of a skew-symmetric monomial of degree $D$ is $d(d - 1)$, where $d = \left[\frac{D}{2}\right]$, i.e. it is bounded from below by a quadratic function in $D$. 
Let us denote this function by $f(D)$. 
Hence, if $w(M) < f(n)$ for $M = x_{k_n}\cdots x_{k_0}$, then $\varphi_2(M) = 0$.

Consider an arbitrary $a\in D_2[\partial]$, $a = a_n\partial^n + \cdots + a_0$. 
Obviously, the weights of the coefficients of $a^N$ are bounded from above by a linear function in $N$ and the degrees of the coefficients of $a^N$ are bounded from below by $N$. 
Hence, there exists such $N$, that $a^N = 0$.
\end{proof}

{\bfseries 4. Some results about the ideal $[x^m]$} 

Using Lemma 1, one can easily deduce an answer to one of Ritt's problems (\cite{Ritt}).

{\bfseries Theorem 3.} \emph{ Let $[x^m]$ be a differential ideal in the algebra of differential polynomials $k\{x\}$ generated by $x^m$ and $x_i$ be the $i$-th derivation of $x$, $i \in \mathbb{Z}_{\geqslant 0}$. Then, the minimal $j$ such that $x_i^j\in [x^m]$ equals $(i + 1)m - i$.}

\begin{proof}

We denote $(i + 1)m - i$ by $q_i$ and claim that $(\varphi_m(\bar{x}_i))^{q_i} = 0$. 
Indeed, $2(m - 1)(i + 1) = 2(q_i - 1)$ basis vectors appear in $\varphi_m(\bar{x}_i)^{q_i}$, but its degree is $2q_i$. 
Since $\varphi_m(\bar{x}_i)^{q_i}$ is a skew-symmetric form, it is equal to zero.

It remains to check that $\varphi_m(\bar{x}_i)^{q_i - 1} \neq 0$. 
Recall that $\bar{x}_i = \sum\limits_{0\leqslant l\leqslant m, 0\leqslant k\leqslant i} C_{i}^k\xi_k^l\land\eta_{i - k}^l $. 
Expanding brackets (without rearranging variables in monomials) in $\varphi_m(\bar{x}_i)^{q_i - 1}$, we obtain a polynomial with positive coefficients. 
It equals $c \prod\limits_{0\leqslant l\leqslant m}\left( \prod\limits_{0\leqslant k\leqslant i}\xi_k^l\land\eta_{i - k}^l \right)$, where $c\in k$ is some constant. 
Clearly, all rearrangings of the variables in the monomials are even permutations. 
Indeed, each rearranging can be performed by moving the whole pair of a form $\xi_k^l\land\eta_{i - k}^l$. Hence, $c$ is positive and nonzero.
\end{proof}

{\bfseries Theorem 4.} (Integral property of the ideal $[x^m]$) \emph{If $f'\in[x^m]$, where $f\in D_m$, then $f\in[x^m]$.}

This theorem was proved earlier by M.V.Kondratieva. We provide an alternative proof.

\begin{proof} 

It is sufficient to prove that there are no constants\footnote{i.e. such elements $C$, that $C' = 0$} in $D_m$. 
Let us prove that there  are no constants in $\Lambda(V_m)$. 
Assume the contrary, then let $C\in\Lambda(V_m)$ be a constant, i.e. $C' = 0$. 
Without loss of generality it can be assumed that $C$ constains $\xi^0_j$ for some $j$. 
Let $j$ be the highest order of a derivation of $\xi^0_0$ in $C$ and $C = \xi^0_j\land a + b$, where $b$ does not contain $\xi^0_j$. 
Then $C' = \xi^0_{j + 1}\land a + c$, where $c$ does not contain $\xi^0_{j + 1}$. 
Hence, $C'\neq 0$. 
This contradiction concludes the proof.

\end{proof}

The author is grateful to A.I.Zobin and Yu.P.Razmyslow for helpful discussions.

\end{document}